\documentclass[12pt,twoside,reqno]{amsart}

\usepackage{color}
\usepackage{graphicx}
\usepackage{amsmath}
\usepackage{amssymb}
\usepackage{amsfonts}
\usepackage{multicol}
\usepackage[ansinew]{inputenc}
\usepackage{amscd}
\usepackage[mathscr]{eucal}
\usepackage{hyperref}
\usepackage{xcolor}
\usepackage{easyReview}

\allowdisplaybreaks

\newcommand{\blob}{\rule[.2ex]{.8ex}{.8ex}}

\newcommand{\R}{\mathbb{R}}

\newcommand{\Z}{\mathbb{Z}}

\newcommand{\SL}{{\rm SL}}

\newcommand{\norm}[1]{\lVert#1\rVert} 
\newcommand{\normtwo}[1]{
	{\left\vert\kern-0.25ex\left\vert\kern-0.25ex\left\vert #1
		\right\vert\kern-0.25ex\right\vert\kern-0.25ex\right\vert} }

\newcommand{\Hcal}{\mathcal{H}}

\newcommand\restr[2]{{
		\left.\kern-\nulldelimiterspace 
		#1 
		\vphantom{\big|} 
		\right|_{#2} 
}}

\theoremstyle{plain}
\newtheorem{theorem}{Theorem}[section]
\newtheorem{proposition}{Proposition}[section]
\newtheorem{corollary}[proposition]{Corollary}
\newtheorem{lemma}[proposition]{Lemma}

\theoremstyle{definition}
\newtheorem{definition}{Definition}[section]

\theoremstyle{definition}
\newtheorem{remark}{Remark}[section]
\newtheorem{example}[theorem]{Example}
\numberwithin{equation}{section}

\newcommand{\Prob}{\mathrm{Prob}}

\title {Moduli of continuity preserved by the Hilbert Transform }

\begin{document}

 \author[P. Duarte]{Pedro Duarte}
 \address{Departamento de Matem\'atica and CEMS.UL\\
 	Faculdade de Ci\^encias\\
 	Universidade de Lisboa\\
 	Portugal
 }
 \email{pmduarte@fc.ul.pt}

 \author[T. Graxinha]{Tom\'e  Graxinha}
 \address{Departamento de Matem\'atica and CEMS.UL\\
 	Faculdade de Ci\^encias\\
 	Universidade de Lisboa\\
 	Portugal}
 \email{graxinhatome@hotmail.com}

 \maketitle

 \begin{abstract}
It has long been known that the Hilbert transform preserves certain  moduli of continuity such as H\"older continuity. In~\cite{GS01}, M. Goldstein and W. Schlag provided a criterion, based on two simple conditions, for a modulus of continuity to be preserved under the Hilbert transform. 
These conditions are not  satisfied by the weak-H\"older and the log-H\"older moduli of continuity.  More recently, A. Avila, Y. Last, M. Shamis and Q. Zhou proved in~\cite{ALSZ21} that for compactly supported $t$-log-H\"older continuous functions ($t>1$) their Hilbert transforms are at least  $(t-1)$-log-H\"older continuous.  A similar result holds for weak-H\"older continuous functions.
In this note we extend these  results to a class of continuous functions satisfying a couple of mild  integrability conditions.
\end{abstract}

 \section{Introduction}

In a recent paper~\cite{dg2026}, the authors constructed a class of random Schrödinger cocycles $A_t$ with non-compact support, in which the failure to satisfy certain moment conditions implies the absence of corresponding moduli of continuity, both for the integrated density of states (IDS) $\rho(t)$ and for the Lyapunov exponent (LE) $L_1(A_t)$ of the random Schrödinger cocycle. This lack of regularity is established first for $\rho(t)$ and subsequently transferred to $L_1(A_t)$ via the Thouless formula:
\[ L_1(A_t) = \int_{-\infty}^\infty \log |s-t| \, d\rho(s) = \pi \cdot \mathcal{H}\rho (t) \]
which can be reinterpreted as the fact that the LE is the Hilbert transform of the IDS, and vice versa (see~\cite{CS83, Dam-survey}).

 In the work~\cite{DG25}, we need to apply Lemma 10.3 from~\cite{GS01}, as well as Corollary 2.3 from~\cite{ALSZ21}, which describe classes of moduli of continuity that are preserved under the Hilbert transform.
 
 The claim in~\cite{GS01} that the Hilbert transform preserves  \textit{weak-H\"older}  moduli of continuity is not entirely correct. This fact was pointed out to one of the authors by M. Shamis in 2022, indicating that it had already been known to some experts for some time. In Proposition 2.2 and Corollary 2.3 of the paper \cite{ALSZ21}, the authors assume that the Hilbert transform is applied to compactly supported functions. This prevents a direct application of their results to Hilbert transforms of Lyapunov exponents of non compactly supported random Schrödinger cocycles.

 The motivations for writing this text, which some experts might find rather monotonous or obvious, are threefold. The first is to clarify a claim made in \cite{GS01} regarding the fact that the Hilbert transform preserves  \textit{weak-H\"older}   moduli of continuity. The second is to extend Corollary 2.3 from \cite{ALSZ21} to functions that do not have compact support. The third is to serve as a supplementary document to be consulted alongside \cite{GS01} and \cite{ALSZ21} in the published version of \cite{DG25}, and potentially, in future work as well.
 
 This document will not be submitted to any journal.

 \bigskip

 \section{A symmetric Schr\"odinger cocycle}
\label{sec 1}

Let $(v_n)_{n\geq 1}$ and $(p_n)_{n\geq 1}$ be  sequences of positive numbers such that  
\begin{itemize}
	\item $\lim_{n\to \infty} v_n=+\infty$;
	\item $\sum_{n=1}^\infty p_n\, \log v_n <\infty$;
	\item $\sum_{n=1}^\infty p_n = 1$.
\end{itemize}
and consider the symmetric measure $\nu\in \Prob(\R)$
\[ \nu = \sum_{n=1}^\infty \frac{p_n}{2}\, \left( \delta_{-v_n} +\delta_{v_n} \right) .\]

Consider the Bernoulli shift
$\sigma:\R^\Z\to \R^\Z$ endowed with the Bernoulli measure
$\nu^\Z$ and the associated (unbounded) Anderson-Bernoulli Schr\"odinger operator $H_\omega:\ell^2(\Z)\to \ell^2(\Z)$,
\[ (H_\omega\psi)_n := -(\psi_{n+1}+\psi_{n-1} )+\omega_n\, \psi_n  \]
where $\omega=(\omega_n)_{n\in \Z}$ is a $\nu^\Z$-typical random  sequence in $\R^\Z$.

Let $A_t:\R^\Z\to \SL_2(\R)$ be the corresponding Schr\"odinger cocycle
\[ A_t(\omega):= \begin{bmatrix}
	\omega_0-t & -1 \\ 1 & 0 
\end{bmatrix}\]
where the parameter $t\in\R$ stands for the particle's energy.
By virtue of the above hypothesis this cocycle has  finite moment
\begin{align}
\int \log\norm{A_t(\omega)} \, d\nu^\Z(\omega) &\leq
\sum_{n=1}^\infty \frac{p_n}{2}\, ( \log v_n  +  \log^+|t| )\nonumber \\
&< C+  \log^+|t| <\infty. \label{finite momemnt const}
\end{align}

\begin{lemma}
	\label{even lemma}
Under the above hypothesis  the first Lyapunov exponent
$t\mapsto L_1(A_t)$ of the Schr\"odinger cocycle is an even function.
\end{lemma}

\begin{proof}
Letting
\[  D:=\begin{bmatrix}
	1 & 0 \\ 0 & -1
\end{bmatrix}\; \text{ and } \;   S(a):=\begin{bmatrix}
	a & -1 \\ 1 & 0
\end{bmatrix}    \]
we have the conjugacy relation
\[D\, S(-a)\, D^{-1} = -S(a) \]
which implies that the products
\[  S(\omega_n)\, \cdots\, S(\omega_2)\, S(\omega_1)\quad \text { and } \quad
S(-\omega_n)\, \cdots\, S(-\omega_2)\, S(-\omega_1) \] 
have the same norm. 
Consider the map $h:\R^\Z\to \R^\Z$, $h(\omega)=-\omega$.  By the previous statement
\[ \norm{ A_t^n(\omega) } = \norm{   A^n_{-t}(h(\omega)) } .\]
Since   $h:\R^\Z\to \R^\Z$  preserves the symmetric Bernoulli measure $\nu^\Z$, i.e., $h_\ast\nu^\Z=\nu^\Z$, taking logarithms, dividing by $n$, integrating w.r.t. $\nu^\Z$ and taking the limit as $n\to \infty$, we get that $L_1(A_t)=L_1(A_{-t})$.
\end{proof}

\bigskip

 \section{Moduli of Continuity}

A function \(\omega \colon [0, 1) \to [0, 
+\infty)\) is called a \textit{modulus of continuity}  (MoC)  provided it is:   
\begin{itemize}
	\item  continuous,
	\item strictly-increasing and 
	\item \(\omega(0) = 0\).
\end{itemize}

Let \((X, d)\) be a metric space.

	A function \(f \colon X \to \mathbb{R}\) is said to have \textit{local MoC}  \(\omega\) if for every $a\in X$, there exist positive constants $r>0$ and $C<\infty$ such that for all $x,y\in X$ with
	$d(x,a),\ d(y,a)<r$,
	\begin{equation*}
		|f(x) - f(y)| \leq C\, \omega(d(x,y)) .
	\end{equation*}

\begin{definition}
	\label{def GS threshold}
Let \( \omega\) be a MoC. We say that $\omega$ is above the GS- threshold \ if    for all $0<\delta<1$

\smallskip

\noindent
\blob\quad $\displaystyle 
\int_0^\delta \frac{\omega(r)}{r} \, dr \lesssim\, \omega(\delta) $

\smallskip

\noindent 
\blob\quad $\displaystyle  \delta\, \int_{\delta}^{\frac{1}{2}}  \frac{\omega (r)}{r^2}   \, dr  \lesssim\,  \omega(\delta)  $
\end{definition}
 
The two conditions above, at the GS threshold, relate to:

\begin{tabbing}
	\quad \= \text{(10.9)}\quad \= $\displaystyle  \sum_{n\geq \ell}   \omega(2^{-n}) \asymp \omega(2^{-\ell})$\\
	\smallskip\\
	\>\text{(10.10)}  \> $\displaystyle   \sum_{n<\ell}    2^n\, \omega(2^{-n}) \asymp 2^\ell\, \omega(2^{-\ell})$
\end{tabbing}	
used in Lemma 10.3 
of the work by Goldstein and Schlag~\cite{GS01}.

Notice that  since $\omega(0)=0$, if $\omega(r)$ is strictly concave then
$r\mapsto \omega(r)/r$ is strictly decreasing and
\begin{align*}
	\int_0^{2^{-\ell}} \frac{\omega(r)}{r}\, dr \leq \sum_{n\geq \ell} \frac{ \omega(2^{-(n+1)})}{2^{-(n+1)}}\, (2^{-n}-2^{-(n+1)}) =   \sum_{n\geq \ell+1}   \omega(2^{-n}) 
\end{align*}
while
\begin{align*}
	\int_0^{2^{-\ell}} \frac{\omega(r)}{r}\, dr \geq \sum_{n\geq \ell} \frac{ \omega(2^{-n})}{2^{-n}}\, (2^{-n}-2^{-(n+1)}) = \frac{1}{2}\, \sum_{n\geq \ell}   \omega(2^{-n}) ,
\end{align*}
so that
\[ 	\int_0^{2^{-\ell}} \frac{\omega(r)}{r}\, dr\asymp  \sum_{n\geq \ell}   \omega(2^{-n}) .\]
Similarly,
\begin{align*}
	\int_{2^{-\ell}}^{\frac{1}{2}} \frac{\omega(r)}{r^2}\, dr \leq \sum_{n=1}^{\ell-1} \frac{ \omega(2^{-(n+1)})}{2^{-2(n+1)}}\, (2^{-n}-2^{-(n+1)}) =   \sum_{n=2}^{\ell}   2^n\, \omega(2^{-n}) 
\end{align*}
while
\begin{align*}
	\int_{2^{-\ell}}^{\frac{1}{2}} \frac{\omega(r)}{r^2}\, dr \geq \sum_{n=1}^{\ell-1} \frac{ \omega(2^{-n})}{2^{-2n}}\, (2^{-n}-2^{-(n+1)}) =   \frac{1}{2}\, \sum_{n=1}^{\ell-1}   2^n\, \omega(2^{-n}) ,
\end{align*}
so that
\[ 	\int_{2^{-\ell}}^{\frac{1}{2}} \frac{\omega(r)}{r^2}\, dr\asymp  \sum_{n< \ell}   2^n\,\omega(2^{-n}) .\]
Notice that the reverse inequalities in conditions (10.9) and (10.10) are obvious.
This proves that the hypothesis at the GS threshold (Definition~\ref{def GS threshold})  are equivalent to the conditions (10.9) and (10.10) used in~\cite[Lemma 10.3]{GS01}. 

\bigskip

One example above the GS-threshold is:

\begin{example}
The \textit{ H\"older} MoC, defined by 
\begin{equation}
	\omega_\alpha(r) :=r^\alpha = \exp\left(-\alpha \log \tfrac{1}{r}\right),
\end{equation}
where $0<\alpha < 1$.
\end{example}

\begin{proof}
	Notice that
	\[  \int_0^\delta \frac{r^\alpha}{r}\, dr = \int_0^\delta r^{\alpha-1}\, dr = \frac{\delta^\alpha}{\alpha} \,  ,  \]
	while
	\[  \delta\, \int_\delta^{\frac{1}{2}} \frac{r^\alpha}{r^2}\, dr = \delta\, \int_\delta^{\frac{1}{2}} r^{\alpha-2}\, dr =  \frac{  ( 2\,\delta)^{\alpha }-2 \delta}{2^\alpha\,(1-\alpha) }\, .  \]
\end{proof}

Another relevant modulus of continuity is

\begin{example}
The \textit{weak-Hölder} MoC, defined by 
\begin{equation}
	\omega_{c,\theta}(r) = \exp\left(- c \left( \log \tfrac{1}{r}\right)^\theta\right),
\end{equation} 
for some $0<c$ and  $0<\theta< 1$.  \\
\end{example}

In~\cite[page 194]{GS01}, Goldstein and Schlag state that the weakly H\"older  MoC lies above the GS threshold, which is not entirely correct.
This fact was pointed out to one of the authors by M. Shamis in 2022, indicating that it had already been known to some experts for some time.
We shall now clarify this issue.

\begin{remark}
Choosing, to simplify the calculations, $c=(\log 2)^{-\frac{1}{2}}$ and $\theta=\frac{1}{2}$,   the weak-H\"older MoC $\omega_{c,\frac{1}{2}}$  does not satisfy condition (10.9) in~\cite{GS01}, because
\[  \sum_{n\geq \ell}   \omega_{c,\frac{1}{2}}(2^{-n}) = \sum_{n\geq \ell} e^{-\sqrt{n}}  \gg   e^{-\sqrt{\ell}} = \omega_{c,\frac{1}{2}}(2^{-\ell}). \]
In fact
\begin{align*}
e^{\sqrt{\ell}}\, 	\sum_{n\geq \ell} e^{-\sqrt{n}} &=
	\sum_{n\geq \ell} e^{\sqrt{\ell} -\sqrt{n}} =\sum_{n\geq \ell} e^{\frac{-(n-\ell)}{\sqrt{\ell} +\sqrt{n}} }\\
	&\geq  \sum_{\sqrt{\ell} \leq k\leq 2\sqrt{ \ell}} e^{\frac{-k}{\sqrt{\ell} +\sqrt{\ell+k}} } \geq \sqrt{\ell}\, e^{-1}
	\to \infty \ \text{ as } \ \ell\to \infty.
\end{align*}
\end{remark}

The family of weak-H\"older MoC nearly satisfies the GS-threshold conditions in the following sense.

\begin{lemma} \label{lemma cond 1 weak holder}
	For $c>0$ and $0< \theta< 1$,  given $0<\delta<1$,
	\begin{align*}
		\int_0^\delta \frac{\omega_{c,\theta}(r)}{r} \, dr \lesssim_{c,\theta}  \, \omega_{\frac{c}{2},\theta}(\delta) .
	\end{align*}  
\end{lemma}

\begin{proof}
We have
\begin{align*}
\int_0^\delta \frac{e^{- c \left( \log \frac{1}{r}\right)^\theta }}{r}\, dr &\leq 
 e^{- \frac{c}{2} \left( \log \frac{1}{\delta}\right)^\theta  }\, \int_0^\delta    \frac{e^{- \frac{c}{2} \left( \log \frac{1}{r}\right)^\theta }}{r}\, dr \\
 &\leq  e^{- \frac{c}{2} \left( \log \frac{1}{\delta}\right)^\theta  }\, \int_{0 }^\infty     e^{- \frac{c}{2} \, u^\theta } \, du \lesssim_{c,\theta} \;  e^{- \frac{c}{2} \left( \log \frac{1}{\delta}\right)^\theta  }.
\end{align*}  
\end{proof}

\begin{lemma} \label{lemma cond 2 weak holder}
		For $c>0$ and $0< \theta< 1$, given $0<\delta< e^{-1}$,
	\begin{equation*}
		\delta\, \int_{\delta}^{\frac{1}{2}}  \frac{\omega_{c,\theta}(r)}{r^2}   \, dr  \lesssim_{c,\theta} \omega_{\frac{c}{2},\theta}(\delta) .
	\end{equation*}  
\end{lemma}

\begin{proof}
	The function $ \omega_{\frac{c}{2},\theta}(r) $ is strictly concave for $0<\delta <e^{-1}$. Therefore, $ \omega_{\frac{c}{2},\theta}(r)/r $ is strictly decreasing as $r$ varies within this interval.
	\begin{align*}
		\int_{\delta}^{\frac{1}{2}}  \frac{\omega_{c,\theta}(r)}{r^2}   \, dr &\leq 
		\frac{\omega_{\frac{c}{2},\theta}(\delta)}{\delta}   \, \int_{\delta}^{\frac{1}{2}}  \frac{e^{- \frac{c}{2} \left( \log \frac{1}{r}\right)^\theta  }}{r}\, dr\\
		& \leq 
		\frac{\omega_{\frac{c}{2},\theta}(\delta)}{\delta}   \,\int_{\log 2}^{\log(\frac{1}{\delta})} e^{-\frac{c}{2}\, u^\theta}\, du \\
		&\leq  \frac{\omega_{\frac{c}{2},\theta}(\delta)}{\delta}\,\int_{0}^{\infty} e^{-\frac{c}{2}\, u^\theta}\, du  
		\lesssim_{c,\theta}\,  \frac{\omega_{\frac{c}{2},\theta}(\delta)}{\delta},
	\end{align*}
	which implies the claimed inequality.
\end{proof}

\bigskip

Consider finally the family of MoC $\omega_\gamma:[0,1)\to [0,\infty)$,
\[ \omega_\gamma(r)=\frac{1}{\left(\log(\frac{1}{ r})\right)^\gamma}\qquad (\gamma>0) \]

We say that a function $f:\R\to \R$ is (locally) $\gamma$-log-H\"older continuous,  if it has \textit{local modulus of continuity} \(\omega_\gamma\).

\begin{lemma} \label{lemma cond 1 omega gamma}
	For $\gamma >1$ and $0<\delta<1$, 
	\begin{equation*}
		\int_0^\delta \frac{\omega_\gamma(r)}{r} \, dr 	= \frac{1}{\gamma-1}\,\omega_{\gamma-1}(\delta) .
	\end{equation*}  
\end{lemma}

\begin{proof}
	A simple calculation gives
	\[ 	\int_0^\delta \frac{\omega_\gamma(r)}{r} \, dr   = \frac{1}{\gamma-1}\,  \frac{1}{(\log (1/\delta ))^{\gamma-1}}
	= \frac{1}{\gamma-1}\,\omega_{\gamma-1}(\delta) . \]
\end{proof}

\begin{lemma} \label{lemma cond 2 omega gamma}
	Given $\gamma >1$  and $0<\delta<e^{-(1+\gamma)}$,
	\begin{equation*}
		\delta\, \int_{\delta}^{\frac{1}{2}}  \frac{\omega_\gamma(r)}{r^2}   \, dr  \lesssim_\gamma  \omega_{\gamma-1}(\delta) .
	\end{equation*}  
\end{lemma}

\begin{proof}
The function $ \omega_{\gamma}(r) $ is strictly concave for $0<\delta <e^{-1-\gamma}$. 
Therefore, $ \omega_{\gamma}(r)/r $ is strictly decreasing as $r$ varies within this interval and
	\begin{align*}
		\int_{\delta}^{\frac{1}{2}}  \frac{\omega_\gamma(r)}{r^2}   \, dr &\leq 
		\frac{\omega_\gamma(\delta)}{\delta}   \, \int_{\delta}^{\frac{1}{2}}  \frac{1}{r}\, dr\leq 
		\frac{\omega_\gamma(\delta)}{\delta}   \,\log\frac{1}{\delta} =  \frac{\omega_{\gamma-1}(\delta)}{\delta},
	\end{align*}
	which implies the claimed inequality.
\end{proof}

\bigskip

\section{MoC preserved by the Hilbert transform}

We define a simple criterion for the Hilbert transform of a continuous function to be well-defined as a Cauchy principal value.
 
 \begin{theorem}
 	\label{thm 1}
Let  $f:\R\to \R$  be a continuous function with $\omega$  as a local MoC such that
\begin{equation} 
	\label{eq MoC condition}
	\int_0^{\frac{1}{2}} \frac{\omega(r)}{r} \, dr <\infty. 
\end{equation}
and assume that
\begin{equation}
	\label{eq 2 integrability conditions}
	\int_{-\infty}^\infty \frac{|f(t)|}{1+t^2}\, dt<\infty \ \text{ and } \
	\int_1^\infty \frac{|f(-t)-f(t)|}{t}\, dt <\infty 
\end{equation}
Then  the Hilbert transform $\rho(s)=(\Hcal f)(s)$
is well defined.
\end{theorem}

\begin{proof}
	Since  
	\[     \mbox{p.v.} \int_{-\infty}^{\infty}\frac{f(x)}{x-t}\, dt = 0, \]
	the Hilbert transform is given by
	\[  (\Hcal f)(x)= \mbox{p.v.} \int_{-\infty}^\infty \frac{f(t)}{x-t}\, dt = \mbox{p.v.} \int_{-\infty}^\infty \frac{f(t)-f(x)}{x-t}\, dt ,\]
	which can be decomposed as follows into a local part (near $x$) and a global part (far from $x$):
	\[ (\Hcal f)(x)= 	\int_{|x-t|\leq \frac{1}{2}} \frac{f(t)-f(x)}{x-t}\, dt \ +	\ \mbox{p.v.}  \int_{|x-t|> \frac{1}{2} } \frac{f(t)-f(x)}{x-t}\, dt .\]

	The local part is  bounded as follows
	\begin{align*}
		\left\vert 	\int_{|x-t|\leq \frac{1}{2} } \frac{f(t)-f(x)}{x-t}\, dt \right\vert &\leq \int_{|x-t|\leq \frac{1}{2}} \frac{|f(t)-f(x)|}{|x-t|}\, dt\\
		&\leq C \, \int_{|x-t|\leq \frac{1}{2}} \frac{\omega(|x-t|)}{|x-t|}\, dt\\
		&\leq 2\, C \,  \int_0^{\frac{1}{2}} \frac{\omega(r)}{r} \, dr \ <\infty .
	\end{align*}
	For the global part, assume first that $f(t)$ is even, i.e., $f(-t)=f(t)$.
	We have in this case
	\begin{align*}
		& \left\vert 		\mbox{p.v.} \int_{|x-t|> \frac{1}{2}} \frac{f(t)-f(x)}{x-t}\, dt \right\vert  = \left\vert 	\int_{|x-t|> \frac{1}{2}} \frac{f(t)}{x-t}\, dt \right\vert \\
		&\qquad \leq \left\vert  \int_{x+\frac{1}{2}}^\infty  \frac{f(t)}{x-t}\, dt 
		+ \int_{-\infty}^{x-\frac{1}{2}}  \frac{f(t)}{x-t}\, dt   \right\vert\\
		&\qquad = \left\vert  \int_{x+\frac{1}{2}}^\infty  \frac{f(t)}{x-t}\, dt 
		+ \int_{\frac{1}{2}-x}^{\infty}  \frac{f(-t)}{x+t}\, dt   \right\vert\\
		&\qquad \leq    \int_{ \frac{1}{2}-|x|}^{\frac{1}{2}+|x|}  \frac{|f(\mp t)|}{|x\pm t|}\, dt  
		+ \left\vert \int_{\frac{1}{2}+|x|}^{\infty}  f(t)\, \left( \frac{1}{x-t} +\frac{1}{x+t}  \right)\, dt   \right\vert\\ 
		&\qquad \leq 4\, |x| \, M_x +
		\int_{\frac{1}{2}+|x|}^{\infty}   \frac{2\, |x|\, |f(t)|}{t^2-x^2}\, dt\\
		&\qquad \leq 2\, |x| \, \left\{  2\,M_x +
		L_x\, \int_{\frac{1}{2}}^{\infty}   \frac{ |f(t)|}{t^2+1}\, dt  \right\} <\infty ,
	\end{align*}
	where $M_x =  \norm{ f\vert_{[\frac{1}{2}-|x|,  \frac{1}{2}+|x|]}}_\infty$ and $L_x =  \sup_{t\geq \frac{1}{2}+|x|} \frac{t^2+1}{t^2- x^2}\leq 5+|x|$.
	This concludes the proof that $\rho(s)=\Hcal f)(s)$ is well-defined.
 
	In general, if $f$ is not an even function, but
	\[ \int_{1}^\infty  \frac{|f(-t)-f(t)|}{t}\, dt <\infty,\]
	it suffices to add (in the third step of the calculation above) the following constant to the upper bound obtained previously:
	\[ \int_{\frac{1}{2}-x}^\infty  \frac{|f(-t)+f(t)|}{x+t}\, dt \lesssim_x  \int_{1/2}^\infty  \frac{|f(-t)-f(t)|}{t}\, dt <\infty .\]
\end{proof}

The following theorem extends Corollary 2.3 from \cite{ALSZ21} and clarifies the assertion that the  weak-H\"older MoC  is preserved by the Hilbert transform, as a consequence of Lemma 10.3 from \cite{GS01}.
\begin{theorem}
	\label{thm 2}
	Let  $f:\R\to \R$  be a continuous function   such that
	\begin{align*}
		\int_{-\infty}^\infty \frac{|f(t)|}{1+t^2}\, dt<\infty \ \text{ and } \
		\int_1^\infty \frac{|f(-t)-f(t)|}{t}\, dt <\infty 
	\end{align*}
 
	\begin{enumerate}
		\item[(a)] If $f(t)$ is locally  $\gamma$-log-H\"older continuous with $\gamma>1$  then $\rho(s)$
		is  locally $(\gamma-1)$-log-H\"older continuous;
		\item[(b)] If $f(t)$ is locally  $(c,\theta)$-weak-H\"older continuous with $0<c\leq 1$ and $0<\theta\leq 1$  then   $\rho(s)$
		is   locally  $(\frac{c}{2},\theta)$-weak-H\"older continuous.
	\end{enumerate}
	
\end{theorem}

\begin{proof}
 Both $\gamma$-H\"older and $(c,\theta)$-weak H\"older moduli of continuity satisfy~\eqref{eq MoC condition}, which allows us to apply Theorem~\ref{thm 1} to conclude that $\rho(s)= (\Hcal f)(s)$ is well defined.

\smallskip

\noindent
{\bf Proof of (a).} \ Consider   $x<y$ in $\R$ such that $|x-y|=\delta <\frac{1}{4}$. Then
\begin{align*}
&\left\vert (\Hcal f)(x) -  (\Hcal f)(y) \right\vert   = 
\left\vert \mbox{p.v.} \int_{-\infty}^\infty \frac{f(t)}{x-t} -  \frac{f(t)}{y-t}\, dt   \right\vert \\
& = 
\left\vert \mbox{p.v.} \int_{-\infty}^\infty \frac{f(t)-f(x)}{x-t} -  \frac{f(t)-f(y)}{y-t}\, dt   \right\vert \\
&\qquad \leq   
\underbrace{  \left\vert \mbox{p.v.} \int_{|x-t|<2\delta} \frac{f(t)-f(x)}{x-t} \right\vert +   \left\vert \mbox{p.v.} \int_{|x-t|<2\delta} \frac{f(t)-f(y)}{y-t}\, dt   \right\vert }_{(I)}  + \\
&\qquad\qquad  \underbrace{ \left\vert \mbox{p.v.} \int_{|x-t|> 2\delta} \frac{f(t)-f(x)}{x-t} -  \frac{f(t)-f(y)}{y-t}\, dt   \right\vert }_{(II)}
\end{align*}
We first deal with the near part (I).
Using Lemma~\ref{lemma cond 1 omega gamma} we have
\begin{align*}
 &\left\vert \mbox{p.v.} \int_{|x-t|<2\delta} \frac{f(t) -f(x)}{x-t}  \right\vert  = \left\vert  \mbox{p.v.} \int_{x-2\delta}^{x+2\delta}  \frac{f(t)-f(x)}{x-t}\, dt  \right\vert   \\
 &\qquad \qquad \leq     \int_{x-2\delta}^{x+2\delta}  \frac{|f(t)-f(x)|}{|x-t|}\, dt  \lesssim 2\, \int_{0}^{2\delta}  \frac{\omega_\gamma(r)}{r}\, dr  \\
 &\qquad \qquad 
 =  \frac{2}{\gamma-1}\,\omega_{\gamma-1}(2\delta) \lesssim_\gamma \, \omega_{\gamma-1}(\delta) .
\end{align*}
Similarly,
\begin{align*}
	&\left\vert \mbox{p.v.} \int_{|x-t|<2\delta} \frac{f(t) -f(y)}{y-t}  \right\vert  = \left\vert  \mbox{p.v.} \int_{y-3\delta}^{y+\delta}  \frac{f(t)-f(y)}{y-t}\, dt  \right\vert   \\
	&\qquad\qquad  \leq     \int_{y-3\delta}^{y+\delta}  \frac{|f(t)-f(y)|}{|y-t|}\, dt     =   \int_{-\delta }^{3\delta}  \frac{\omega_\gamma(|r|)}{|r|}\, dr \\
	&\qquad\qquad  = \int_{0}^{\delta}  \frac{\omega_\gamma(r)}{r}\, dr + \int_{0}^{3\delta}  \frac{\omega_\gamma(r)}{r}\, dr \\
	&\qquad\qquad =  \frac{1}{\gamma-1}\,\left( \omega_{\gamma-1}(\delta) +  \omega_{\gamma-1}(3\delta) \right) \lesssim_\gamma   \,\omega_{\gamma-1}(\delta) .
\end{align*}
Hence, combining these two estimates we get \
 $(I)	\;  \lesssim_\gamma   \,\omega_{\gamma-1}(\delta)$.

\bigskip

For part (II), 
using Lemma~\ref{lemma cond 2 omega gamma} and that for $y$ close enough to $x$, say $\delta =|x-y|<\frac{1}{4}$, we have
\[  \left\vert\frac{1}{x-t} -  \frac{1}{y-t}\right\vert\leq  2\, \frac{|x-y|}{(x-t)^2} ,\]
and noticing that  $\{t\in \R\colon |y-t|> 3\delta \}\subset \{t\in \R\colon |x-t|> 2\delta \}$, with difference set $[y+\delta, y+3\delta[$,    and that
\[  \mbox{p.v.} \int_{|y-t|> 3\delta} \frac{ f(x)}{y-t} \, dt  \ = \ \mbox{p.v.} \int_{|y-t|> 3\delta} \frac{ f(y)}{y-t} \, dt \ = 0 , \]
we get
\begin{align*}
(II)&\leq \left\vert \mbox{p.v.} \int_{|x-t|> 2\delta} \frac{f(t)-f(x)}{x-t} -  \frac{f(t)-f(y)}{y-t}\, dt   \right\vert  \\
&\leq \left\vert \mbox{p.v.} \int_{|y-t|> 3\delta} \frac{f(t)-f(x)}{x-t} -  \frac{f(t)-f(y)}{y-t}\, dt   \right\vert \\
&\qquad\qquad +
\left\vert  \int_{y+\delta}^{y+3\delta}   \frac{f(t)-f(y)}{y-t}\, dt   \right\vert   \\
&\leq \left\vert \mbox{p.v.} \int_{|y-t|> 3\delta} \frac{f(t)-f(x)}{x-t} -  \frac{f(t)-f(x)}{y-t}\, dt   \right\vert \\
&\qquad\qquad +
  \int_{y+\delta}^{y+3\delta}   \frac{|f(t)-f(y)|}{|y-t|}\, dt       \\
&\leq  \int_{|y-t|> 3\delta} |f(t)-f(x)|\, \left\vert\frac{1}{x-t} -  \frac{1}{y-t}\right\vert \, dt   \, +   \,\int_0^{3\delta}  \frac{\omega_\gamma(r)}{r} \, dr   \\
&\lesssim \, \int_{|y-t|> 3\delta} |f(t)-f(x)|\,  \frac{|x-y|}{(x-t)^2}   \, dt  \, + \,   \omega_{\gamma-1}(3\delta)  \\
&\leq  \delta\, \left\{  \int_{[-\frac{1}{2},\frac{1}{2}] \setminus [-2\delta, 4\delta] } \frac{\omega_\gamma(|r|)}{r^2}   \, dr   + 2 \int_{\frac{1}{2}} ^\infty \frac{|f(|t|)|+|f(x)|}{t^2} \, dt \right\}  \\
& \hspace{8cm}  +\,  \omega_{\gamma-1}(\delta)  \\
 &\lesssim_\gamma 2\,  \omega_{\gamma-1}(4 \delta) + 2\, \delta \, C +   \omega_{\gamma-1}( \delta) \,  \lesssim_\gamma \,  \omega_{\gamma-1}(\delta).
\end{align*}
This proves that $\Hcal f$ is locally $(\gamma-1)$-log-H\"older continuous.

\smallskip

\noindent
{\bf Proof of (b).} \  The proof is entirely analogous, using Lemma~\ref{lemma cond 1 weak holder} instead of Lemma~\ref{lemma cond 1 omega gamma}, and  Lemma~\ref{lemma cond 2 weak holder} instead of Lemma~\ref{lemma cond 2 omega gamma}.
\end{proof}

\bigskip

\begin{corollary} Consider the first Lyapunov exponent
	$f(t)= L_1(A_t)$ of the model introduced in Section~\ref{sec 1}.
	
\begin{enumerate}
	\item[(a)] 	If  $f$  is locally  $\gamma$-log-H\"older continuous with $\gamma>1$ then the Hilbert transform $\rho(s)=(\Hcal f)(s)$
	is well defined and is, locally,   $(\gamma-1)$-log-H\"older continuous.
		
	\item[(a)] 	If  $f$  is locally  $(c,\theta)$-weak-H\"older continuous with $0<c,\theta\leq1$ then the Hilbert transform $\rho(s)=(\Hcal f)(s)$
is well defined and is, locally,   $(\frac{}{2},\theta)$-weak-H\"older continuous.
\end{enumerate}
	
\end{corollary}

\begin{proof}
By Lemma~\ref{even lemma}, the function $f(t)=L_1(A_t)$ is even,
 which makes the second integrability condition~\eqref{eq 2 integrability conditions} trivial.
 
The moment condition~\eqref{finite momemnt const} implies that
\[  \int_{-\infty}^\infty \frac{|f(t)|}{1+t^2}\, dt \leq 	\int_{-\infty}^\infty \frac{C+\log^+ |t|}{1+t^2}\, dt <\infty .\]
Therefore, the hypotheses of Theorem~\ref{thm 2} are satisfied and the claimed conclusions follow.
\end{proof}

\bigskip

\subsection*{AI Declaration:}

The authors used ChatGPT 5.5 to analyze the proof of Lemma 10.3 in~\cite{GS01}. The use of integral versions of conditions (10.9) and (10.10) from~\cite{GS01}, as well as the more elementary approach, were suggested by this AI model. The authors are responsible for the lengthy calculations, that is, they were not performed by the AI model.


\begin{thebibliography}{1}
	
	\bibitem{ALSZ21}
	Artur Avila, Yoram Last, Mira Shamis, and Qi~Zhou, \emph{On the abominable
		properties of the almost {M}athieu operator with well-approximated
		frequencies}, Duke Math. J. \textbf{173} (2024), no.~4, 603--672.
	\MR{4734551}
	
	\bibitem{CS83}
	Walter Craig and Barry Simon, \emph{Log {H}\"older continuity of the integrated
		density of states for stochastic {J}acobi matrices}, Comm. Math. Phys.
	\textbf{90} (1983), no.~2, 207--218. \MR{714434}
	
	\bibitem{Dam-survey}
	David Damanik, \emph{Schr\"odinger operators with dynamically defined
		potentials}, Ergodic Theory Dynam. Systems \textbf{37} (2017), no.~6,
	1681--1764. \MR{3681983}
	
	\bibitem{DG25}
	Pedro Duarte and Tom\'{e} Graxinha, \emph{H\"older continuity of {L}yapunov
		exponents for non-invertible and non-compact random cocycles}, 2025, preprint
	arXiv:2506.04124 [math.DS].
	
	\bibitem{dg2026}
	\bysame, \emph{Obstructions to the regularity of the lyapunov exponents for
		non-compact random schr\"odinger cocycles}, 2026.
	
	\bibitem{GS01}
	M.~Goldstein and W.~Schlag, \emph{H{\"o}lder continuity of the integrated
		density of states for quasi-periodic {S}chr{\"o}dinger equations and averages
		of shifts of subharmonic functions}, Annals of Math. \textbf{154} (2001),
	155--203.
	
\end{thebibliography}

\providecommand{\bysame}{\leavevmode\hbox to3em{\hrulefill}\thinspace}
\providecommand{\MR}{\relax\ifhmode\unskip\space\fi MR }
\providecommand{\MRhref}[2]{%
	\href{http://www.ams.org/mathscinet-getitem?mr=#1}{#2}
}
\providecommand{\href}[2]{#2}

\end{document}